\documentclass[12pt, reqno]{article}

\usepackage[usenames]{color}
\usepackage{amssymb}
\usepackage{amsmath}
\usepackage{amsthm}
\usepackage{amsfonts}
\usepackage{amscd}
\usepackage{graphicx}
\usepackage{mathtools}

\usepackage[colorlinks=true,
linkcolor=webgreen,
filecolor=webbrown,
citecolor=webgreen]{hyperref}

\definecolor{webgreen}{rgb}{0,.5,0}
\definecolor{webbrown}{rgb}{.6,0,0}

\usepackage{color}
\usepackage{fullpage}
\usepackage{float}

\usepackage{graphics}
\usepackage{latexsym}
\usepackage{epsf}

\usepackage[colorinlistoftodos]{todonotes}

\newcommand{\seqnum}[1]{\href{https://oeis.org/#1}{\rm \underline{#1}}}

\usepackage{multicol}
\newcommand{\abs}[1]{\lvert#1\rvert}

\begin{document}

\theoremstyle{plain}
\newtheorem{theorem}{Theorem}
\newtheorem{corollary}[theorem]{Corollary}
\newtheorem{lemma}[theorem]{Lemma}
\newtheorem{proposition}[theorem]{Proposition}

\newtheorem*{boundon-quatupperdensity}{Theorem~\ref{mainthm-quatupperlowerbd}}
\newtheorem*{densityRankin}{Theorem~\ref{Theorem:Rankindensity}}

\theoremstyle{definition}
\newtheorem{definition}[theorem]{Definition}
\newtheorem{example}[theorem]{Example}
\newtheorem{conjecture}[theorem]{Conjecture}

\newtheorem{question}[theorem]{Question}

\theoremstyle{remark}
\newtheorem{remark}[theorem]{Remark}

\begin{center}
\vskip 1cm{\Large \textbf{Avoiding 3-Term Geometric Progressions in \\
Hurwitz Quaternions}}
\end{center}

\vskip .1in

\begin{multicols}{2}
\begin{center}
Megumi Asada \\
Graduate School of Education \\
Rutgers University \\
New Brunswick, NJ 08901\\
USA\\
\href{mailto:megumi.asada@rutgers.edu}{megumi.asada@rutgers.edu}\\
\end{center}
\columnbreak
\begin{center}
Bruce Fang \\
Department of Mathematics\\
and Statistics \\
Williams College \\
Williamstown, MA 01267 \\
USA \\
\href{mailto:bf8@williams.edu}
{bf8@williams.edu}
\end{center}
\end{multicols}

\begin{multicols}{2}
\begin{center}
Eva Fourakis \\
Department of Mathematics\\
and Statistics \\
Williams College \\
Williamstown, MA 01267 \\
USA \\
\href{mailto:erf1@williams.edu}{erf1@williams.edu}
\end{center}
\columnbreak
\begin{center}
Sarah Manski \\
Center for Statistical Training\\
and Consulting \\
Michigan State University \\
East Lansing, MI 48824 \\
USA \\
\href{mailto:mailto:manskisa@msu.edu}{manskisa@msu.edu}
\end{center}
\end{multicols}

\begin{multicols}{2}
\begin{center}
Nathan McNew \\
Department of Mathematics \\
Towson University \\
Towson, MD 21252 \\
USA \\
\href{mailto:nmcnew@towson.edu}{nmcnew@towson.edu}
\end{center}

\begin{center}
\columnbreak
Steven J. Miller \\
Department of Mathematics\\
and Statistics \\
Williams College \\
Williamstown, MA 01267 \\
USA \\
\href{mailto:sjm1@williams.edu}{sjm1@williams.edu} \\
\href{Steven.Miller.MC.96@aya.yale.edu}{Steven.Miller.MC.96@aya.yale.edu}
\end{center}
\end{multicols}

\begin{multicols}{2}
\begin{center}
Gwyneth Moreland \\
Department of Mathematics, Statistics, and Computer Science \\
University of Illinois Chicago\\
Chicago, IL 60607\\
USA\\
\href{mailto:gwynm@uic.edu}{gwynm@uic.edu}
\end{center}
\columnbreak
\begin{center}
Ajmain Yamin \\
Department of Mathematics \\
CUNY Graduate Center \\
New York, NY 10016 \\
USA\\
\href{mailto:ayamin@gradcenter.cuny.edu}{ayamin@gradcenter.cuny.edu}
\end{center}
\end{multicols}

\vskip .1in

\begin{center}
Sindy Xin Zhang \\
Department of Mathematics \\
Tufts University \\
Medford, MA 02155 \\
USA\\
\href{mailto:sindy.zhang@tufts.edu}
{sindy.zhang@tufts.edu}
\end{center}

\begin{abstract}

Several recent papers have considered the problem of how large a subset of integers can be without containing any 3-term geometric progressions. This problem has also recently been generalized to rings of integers in quadratic number fields and polynomial rings over finite fields. We study the analogous problem in the Hurwitz quaternion order to see how non-commutativity affects the problem.  We compute an exact formula for the density of a 3-term geometric-progression-free set of Hurwitz quaternions arising from a greedy algorithm and derive upper and lower bounds for the supremum of upper densities of 3-term geometric-progression-free sets of Hurwitz quaternions.
\end{abstract}

\section{Introduction}\label{sec:one}
In 1961, Rankin \cite{Ran} introduced the problem of finding large sets of positive integers which avoid $3$-term geometric progressions.  An obvious example of such a set is the set of square-free positive integers which has asymptotic density $6/\pi^2 \approx 0.607927$.  Rankin constructed a higher density geometric-progression-free set, which we denote as $G_3^*(\mathbb{N}_+)$, and obtained an exact closed-form expression for its density in terms of the Riemann $\zeta$-function
$$d(G_3^*(\mathbb{N}_+)) \ = \ \frac{1}{\zeta(2)} \prod_{n=1}^\infty \frac{\zeta(3^n)}{\zeta(2\cdot 3^n)} \ \approx \ 0.719745.$$
It was shown by Brown and Gordon \cite{BG} that $G_3^*(\mathbb{N}_+)$ is actually generated by a greedy algorithm; it is formed from the singleton $\{1\}$ by greedily adjoining to it integers of increasing magnitude so long as the enlarged set still avoids 3-term geometric progressions.  As such, we refer to $G_3^*(\mathbb{N}_+)$ as \emph{Rankin's greedy set}.  Listing the elements of $G_3^*(\mathbb{N}_+)$ in increasing order yields the sequence $(1, 2, 3, 5, 6, 7,\ldots)$, which is OEIS \seqnum{A000452} \cite{OEIS}. 

By modifying Rankin's greedy set, McNew \cite{McN} constructed a set of positive integers avoiding 3-term geometric progressions with asymptotic density greater than $d(G_3^*(\mathbb{N}_+))$.  It remains an open problem to determine the supremum of asymptotic densities of such sets.  

Variants of this problem have been studied by many authors including Riddell \cite{Ri},  Brown and Gordon \cite{BG},  Beiglb\"{o}ck, Bergelson, Hindman, 
and Strauss \cite{BBHS}, Nathanson and O'Bryant \cite{NO}, and McNew \cite{McN}.  More recently, Best,
Huan, McNew, Miller, Powell, Tor,  and Weinstein \cite{BHMMPTW} studied analogous problems in rings of integers of quadratic number fields, while Asada, Fourakis, Manski, McNew, Miller, and Moreland \cite{AFMMMM} tackled the problems in polynomial rings over finite fields.

The purpose of this paper is to investigate this problem in a non-commutative setting, namely the Hurwitz quaternion order $Q_{\textup{Hur}}$ (see Section \ref{sec:two}). We will take a \emph{3-term geometric progression} of Hurwitz quaternions to mean a triple $(a,ar,ar^2)$, where $a, r \in Q_\textup{Hur} \setminus \{ 0 \}$ and $r$ a non-unit.  A set $A \subset Q_\textup{Hur}$ \emph{avoids} 3-term geometric progressions if $A^3$ does not contain any such triples.
We are interested in how large such a 3-term geometric-progression-free set can be.

Throughout the paper, when $A \subset \mathbb{N}_+$ is any set of positive integers, we will denote by $S(A) = \{q \in Q_{\textup{Hur}} \mid \textup{Nm}(q) \in A\}$ the set of Hurwitz quaternions whose norm is contained in the set $A$.  We will also abuse this notation slightly, writing $S(N)$ for an integer $N$ to mean $S(N) = S([1,N])$, the set of all non-zero Hurwitz quaternions whose norm is at most $N$. 

 It is easy to see that if $G$ is any set of positive integers which avoid $3$-term geometric progressions, then the set $S(G)$ of Hurwitz quaternions is also 3-term geometric-progression-free, in the sense above.  Our first result is an exact formula for the density (see Section \ref{Section: Rankin}) of $S(G)$ in $Q_\textup{Hur}$ when $G$ is Rankin's greedy set $G_3^*(\mathbb{N}_+)$.

\begin{densityRankin}
The density of the set $S(G_3^*(\mathbb{N}_+))$ of Hurwitz quaternions with norm in Rankin's greedy set $G_3^*(\mathbb{N}_+)$ is 
$$d(S(G_3^*(\mathbb{N}_+))) \ = \ \frac{f(2^2)}{\zeta(2)}  \prod_{\substack{p >2 \\  \emph{prime}}}\left(f(p)-\frac{f(p^2)}{p}\right) \ \approx \ 0.782643,$$ where $$f(x) \ = \ \prod_{i=0}^\infty \left(1+\frac{1}{x^{3^i}}\right).$$

\end{densityRankin}
This is somewhat higher than the density of Rankin's greedy set $G_3^*(\mathbb{N}_+)$ in $\mathbb{N}_+$, which is $\approx 0.719745$.

Our second result provides upper and lower bounds on the supremum of upper densities (see Section \ref{sec:three}) of sets of Hurwitz quaternions free of 3-term geometric progressions.
\begin{boundon-quatupperdensity}
The supremum $\overline{m}_{\textup{Hur}}$ of upper densities of sets of Hurwitz quaternions that avoid 3-term geometric progressions satisfies the following bounds:
$$ 0.946589 \ \approx \ \frac{17665627}{18662400} \ \le\ \overline{m}_{\textup{Hur}} \ \le\ \frac{20}{21} \ \approx \ 0.952381.$$
\end{boundon-quatupperdensity}

In Section \ref{sec:two}, we establish some asymptotic formulas for certain counts of Hurwitz quaternions which are useful for computing densities.  In Section \ref{Section: Rankin}, we prove Theorem \ref{Theorem:Rankindensity}.  The approach is to first prove Theorem \ref{thm:eulerProdForQuatDens}, which is a general Euler product formula for densities of certain sets of the form $S(G)$, and then apply it to the special case when $G$ is Rankin's greedy set $G_3^*(\mathbb{N}_+)$.  In Section \ref{sec:three}, we prove Theorem \ref{mainthm-quatupperlowerbd}.  The approach is to adapt and modify the arguments of McNew \cite{McN}, which are used to produce the best bounds on upper densities of geometric-progression-free sets of integers, to work in the quaternionic setting.  The remainder of the paper is dedicated to the quaternion greedy set $G_3^*(Q_\textup{Hur}\setminus\{0\})$.  A major difference between the rational integer setting and the Hurwitz quaternion setting is the fact that the natural analogue of Rankin's greedy set for $Q_\textup{Hur}$ is much more difficult to analyze, due to the non-commutativity of the quaternions, which was the primary motivation for this work---to investigate what happens in this more challenging situation.  In Section \ref{Section: Greedy Set} we explicate these difficulties and in Section \ref{section:FutureWork} we list further open questions and ways to generalize this work.

\section{Combinatorics of Hurwitz quaternions}\label{sec:two}

The \textit{Hamilton quaternions}
$\mathbb{H} \coloneqq\{a+bi+cj+dk \mid a,b,c,d\in \mathbb{R} \}$ constitute the non-commutative algebra over the reals generated by units $i$, $j$, and $k$ satisfying 
\begin{gather*}
i^2 \ = \  j^2 \ = \  k^2 \ = \  ijk \ = \  -1,\\
ij \ = \ -ji \ = \ k, \hspace{0.5cm} jk \ = \ -kj \ = \ i, \hspace{0.5cm} ki \ = \ -ik \ = \ j.
\end{gather*}
The \textit{norm of a quaternion} $Q = a + bi + cj + dk$ is $\textup{Nm}(Q) \coloneqq a^2+b^2+c^2+d^2$. A quaternion $a + bi + cj + dk$ belongs to the \emph{Hurwitz order} $Q_\text{Hur}$, or is a \emph{Hurwitz quaternion}, if $a,b,c,d$ are all in $\mathbb{Z}$ or all in $\mathbb{Z} + \frac{1}{2}$. The norm of a Hurwitz quaternion is always a nonnegative integer.

In this section, we compute asymptotic formulas for the number of Hurwitz quaternions up to a certain norm excluding 0, and the subset of those whose norm is coprime to a given number. These formulas will be used in the proofs of our main results Theorem \ref{Theorem:Rankindensity} and Theorem \ref{mainthm-quatupperlowerbd}.  For a more in-depth discussion of the Hurwitz quaternion order, see Conway and Smith \cite[Ch.\ 5]{ConSm} or Voight \cite[Ch.\ 11]{V}.

We begin with a classical result that follows from Jacobi's four-square theorem.
\begin{lemma}\label{Lemma:NumquatsfixedN} The number of Hurwitz quaternions of norm $N$ is $24\cdot \sigma_{\emph{odd}}(N)$, where $\sigma_{\emph{odd}}$ is the sum-of-odd-divisors function $$\sigma_{\emph{odd}}(N) \ \coloneqq \ \sum_{2\nmid d \mid N} d.$$
\end{lemma}
For a complete proof of this fact see, for example, the master's thesis of Negrini \cite{N}.  The sequence of numbers $( 24 \cdot \sigma_\text{odd}(N) )_{N>0} = (24, 24, 96, 24, 144, \ldots )$ is OEIS \seqnum{A004011} \cite{OEIS_2}.

Using this lemma, we prove the following asymptotic result.
\begin{lemma}\label{Lemma:NumquatsleN} \leavevmode
\begin{enumerate}
\item The number of Hurwitz quaternions with norm less than or equal to $N$ is
$$\sum_{n\leq N}24\cdot \sigma_{\emph{odd}}(n) \ = \  \pi^2 N^2 + O(N\log N).$$
\item For a fixed integer $M$, the number of Hurwitz quaternions with norm at most $N$ and with norm coprime to $M$ is
$$\sum_{\substack{n\leq N\\ \gcd(n,M) = 1}} 24\cdot \sigma_{\emph{odd}}(n) \ = \ \pi^2N^2 \cdot \frac{\varphi(2M)}{M}\prod_{p\mid M}(1-p^{-2}) + O(N\log N).$$
\end{enumerate}
\end{lemma}

\begin{proof}[Proof of (1)]
Making use of the fact that the sum of the first $k$ odd numbers is $k^2$, we compute
\begin{align*}
\sum_{n\leq N} \sigma_{\text{odd}}(n) &\ = \ \sum_{m\leq N} \sum_{\substack{d|m \\ 2\nmid d}} d  \ = \ \sum_{e \leq N} \sum_{\substack{d\leq \left \lfloor \!\! \frac{N}{e}\!\! \right \rfloor \\ 2 \nmid d}} d \tag{write $m=de$} = \ \sum_{e\leq N} \left\lceil \frac{\left \lfloor\frac{N}{e}\right\rfloor}{2}\right\rceil^2  \\
&= \ \frac{1}{4}\sum_{e\leq N}  \left ( \left( \frac{N}{e} \right) ^2 + O\left( \frac{N}{e} \right) \right)  = \ \frac{1}{4}N^2 \sum_{e\leq N} \frac{1}{e^2} + O \left( N\log N \right) \\
&\ = \ \frac{1}{4}N^2  \left(\frac{\pi^2}{6}\right)+ O \left(N\log N + N^2\sum_{e>N}\frac{1}{e^2} \right) \tag{first series converges to $\frac{\pi^2}{6}$}\\
&\ = \  \frac{\pi^2}{24}N^2+ O \left( N\log N\right). \tag*{\qedhere} 
\end{align*}
\end{proof}

\begin{proof}[Proof of (2)] We compute
\begin{align*}
\sum_{\substack{n\leq N \\ \gcd(n,M)=1}}\sigma_{\text{odd}}(n) &\ = \ \sum_{\substack{n\leq N\\\gcd(n,M)=1}}\sum_{2\nmid d\mid n} d \ = \ \sum_{\substack{e\leq N\\ \gcd(e,M)=1}} \sum_{\substack{d\leq \left\lfloor\frac{N}{e}\right\rfloor\\\gcd(d,2M)=1}}d \\ 
&\ = \ \sum_{\substack{e\leq N \\ \gcd(e,M)=1}} \sum_{\substack{r\leq 2M\\ \gcd(r,2M)=1}} \sum_{\substack{d\leq \left\lfloor\frac{N}{e}\right\rfloor\\ d\equiv r \textup{ (mod $2M$)}}} d.
\end{align*}

Since
\begin{align*}
\sum_{\substack{e\leq N \\ \gcd(e,M)=1}} \sum_{\substack{r\leq 2M\\ \gcd(r,2M)=1}} 
\left(\sum_{\substack{d\leq \left\lfloor\frac{N}{e}\right\rfloor\\ d\equiv r \textup{ (mod $2M$)}}} d  \ - \   \sum_{\substack{m\leq \left\lfloor \frac{N}{e}\right\rfloor\\ 2M\mid m}} m\right)   \ = \ O(N\log N),
\end{align*}
we see
\begin{align*}
\sum_{\substack{n\leq N \\ \gcd(n,M)=1}}\sigma_{\text{odd}}(n) 
&\ = \ \sum_{\substack{e\leq N \\ \gcd(e,M)=1}} \sum_{\substack{r\leq 2M\\ \gcd(r,2M)=1}} \sum_{\substack{m\leq \left\lfloor \frac{N}{e}\right\rfloor\\ 2M\mid m}} m  \quad~+~\quad O(N\log N) \\
&\ = \ \varphi(2M)\sum_{\substack{e\leq N \\ \gcd(e,M)=1}} \sum_{n\leq \big\lfloor\frac{\left\lfloor \frac{N}{e} \right\rfloor}{2M} \big\rfloor} 2M \cdot n
\quad ~+~ \quad O(N\log N) \\
&\ = \ M \varphi(2M) \sum_{\substack{e\leq N \\ \gcd(e,M)=1}} \bigg\lfloor\frac{\left\lfloor \frac{N}{e} \right\rfloor}{2M} \bigg\rfloor \left(\bigg\lfloor\frac{\left\lfloor \frac{N}{e} \right\rfloor}{2M} \bigg\rfloor +1\right)\quad~+~ \quad O(N\log N) \\
&\ = \ N^2 \cdot  \frac{\varphi(2M)}{4M} \sum_{\substack{1\leq e\leq \infty\\ \gcd(e,M)=1}} \frac{1}{e^2} \quad~+~ \quad O(N\log N).
\end{align*}
By the inclusion-exclusion principle, expressed with the M\"{o}bius $\mu$ function, we have
\begin{align*}
\sum_{\substack{1\leq e\leq \infty\\ \gcd(e,M)=1}} \frac{1}{e^2} = \   \sum_{d\mid M} \left( \mu(d) \cdot \sum_{n=1}^\infty \frac{1}{(dn)^2}\right)
&\ = \  \zeta(2) \sum_{d\mid M} \frac{\mu(d)}{d^2}   \ = \  \frac{\pi^2}{6} \prod_{p\mid M} (1- p^{-2}),
\end{align*}
and so
\begin{align*}
\sum_{\substack{n\leq N \\ \gcd(n,M)=1}}\sigma_{\text{odd}}(n) &\ =  \ N^2 \cdot \frac{\pi^2}{24} \frac{\varphi(2M)}{M} \prod_{p\mid M} (1- p^{-2}) \quad~+~ \quad O(N\log N).
\qedhere\end{align*}
\end{proof}

\section{Exact densities of sets of Hurwitz quaternions}\label{Section: Rankin}
\begin{definition}
Let $A\subset Q_{\textup{Hur}}$.   The \textit{density of $A$ in $Q_{\textup{Hur}}$}, denoted $d(A)$, is defined to be
\begin{align}\label{eq:densityDef}
d(A) \ \coloneqq \ \lim_{N\rightarrow\infty}\frac{\abs{A\cap S(N)}}{\abs{S(N)}},
\end{align}
provided the limit exists.  Here, as before, $S(N)$ is the set of Hurwitz quaternions of norm at most $N$ (excluding 0).
\end{definition}

\begin{example}
If $A$ denotes the set of Hurwitz quaternions with norm coprime to $M$, then by Lemma \ref{Lemma:NumquatsleN}, we see that the density of $A$ in $Q_{\text{Hur}}$ is $$d(A) \ = \ \frac{\varphi(2M)}{M}\prod_{p\mid M} (1 - p^{-2}).$$
\end{example}

The purpose of this section is to compute an exact formula for the density of $S(G_3^*(\mathbb{N}_+))$ in $Q_{\text{Hur}}$.  Recall that $S(G_3^*(\mathbb{N}_+))$ is the set of Hurwitz quaternions with norm in Rankin's greedy set $G_3^*(\mathbb{N}_+)$.  The set $G_3^*(\mathbb{N}_+)$ is formed by starting with the singleton $\{1\}$ and greedily adjoining to it larger positive integers which do not form 3-term geometric progressions with the previously included elements.  Listing the elements of $G_3^*(\mathbb{N}_+)$ in increasing order yields the sequence $(1,2,3,5,6,7,8,10,11, \ldots)$, which is OEIS \seqnum{A000452} \cite{OEIS}.  Note that 
$$G_3^*(\mathbb{N}_+)=\{n\in\mathbb{N}_{+}\mid \nu_p(n)\in A_3^*(\mathbb{N}_+)\text{ for all primes }p\},$$ where $\nu_p$ denotes the $p$-adic valuation and $A_3^*(\mathbb{N}_+)$ is the set formed by greedily including non-negative integers that do not form a $3$-term arithmetic progression with the previously included elements. Brown and Gordon \cite{BG} showed that $A_3^*(\mathbb{N}_+)$ is also the set of non-negative integers whose ternary expansion does not contain the digit 2.  Listing the elements of $A_3^*(\mathbb{N}_+)$ in increasing order yields the sequence $(0,1,3,4,9,10,12,13,\ldots)$, which is OEIS \seqnum{A005836} \cite{OEIS_4}.

We will compute the density of $S(G_3^*(\mathbb{N}_+))$ in $Q_\text{Hur}$ using the following theorem, which provides an Euler product for densities of certain kinds of sets in $Q_{\text{Hur}}$.  We would like to thank Emma Dinowitz for helpful discussions regarding the proof of the following theorem.

\begin{theorem}\label{thm:eulerProdForQuatDens}
Let $A \subset \mathbb{N}$ and let $G \coloneqq \{n\in \mathbb{N}_{+} \mid \nu_p(n) \in A \textup{ for all primes } p\}$.  Then
$$d(S(G)) \ = \ \frac{1}{\zeta(2)} \left(\sum_{k\in A} \frac{1}{2^{2k}}\right)\prod_{\substack{p>2\\ \emph{prime}}}\left(\sum_{k\in A} \frac{1}{p^k} - \frac{1}{p}\sum_{k\in A} \frac{1}{p^{2k}}\right).$$
\end{theorem}

\begin{proof}
For each prime $p$ and each positive integer $N$, we define $$G_p \ \coloneqq \ \{n\in \mathbb{N}_+ \mid \nu_p(n) \in A\} \quad \text{ and } \quad G_N \ \coloneqq \ \bigcap_{\substack{p\mid N\\ p \text{ prime}}}G_p.$$  For each positive integer $r$, write $N_r$ to denote the $r^{\text{th}}$ \emph{primorial} $N_r \coloneqq p_1p_2\cdots p_r$.  For convenience, write $\Sigma(N)$ to denote $\Sigma(N) \coloneqq \sum_{n\leq N} \sigma_{\text{odd}}(n)$.  We compute
\begin{align*}
d(S(G)) &\ = \ \lim_{N\rightarrow \infty} \frac{1}{\sum_{n\leq N} \sigma_{\text{odd}}(n)} \sum_{n\leq N} \sigma_{\text{odd}}(n) \cdot \mathbf{1}_{G}(n) \\
&\ = \ \lim_{r\rightarrow \infty} \lim_{N\rightarrow \infty} \frac{1}{\Sigma(N)} \sum_{n\leq N} \sigma_{\text{odd}}(n)\cdot \mathbf{1}_{G_{N_r}}(n) \\
&\ = \ \lim_{r\rightarrow \infty} \lim_{N\rightarrow \infty} \frac{1}{\Sigma(N)} \sum_{(a_j)_{j=1}^r \in \mathbb{N}^r}  \sum_{\substack{n\leq N\\ \forall j\leq r,\\ \nu_{p_j}(n) = a_j}}\sigma_{\text{odd}}(n) \cdot \mathbf{1}_{G_{N_r}}(p_1^{a_1} \cdots p_r^{a_r}).
\end{align*}
Since $\sigma_{\text{odd}}$ is multiplicative, we see
\begin{align*}
d(S(G)) &\ = \ \lim_{r\rightarrow \infty} \lim_{N\rightarrow \infty} \frac{1}{\Sigma(N)} \sum_{(a_j)_{j=1}^r \in \mathbb{N}^r}  \sum_{\substack{n\leq N\\ \forall j\leq r,\\ \nu_{p_j}(n) = a_j}}\sigma_{\text{odd}}\left(\frac{n}{p_1^{a_1}
\cdots p_r^{a_r}}\right)  \\
& \quad \quad \quad \quad \quad \quad \quad \quad \quad \quad \quad \quad \quad \cdot \sigma_{\text{odd}}(p_1^{a_1}\cdots p_r^{a_r})\mathbf{1}_{G_{N_r}}(p_1^{a_1}\cdots p_r^{a_r}) \\
&\ = \ \lim_{r\rightarrow \infty} \sum_{(a_j)_{j=1}^r \in \mathbb{N}^r}  \lim_{N\rightarrow \infty} \frac{1}{\Sigma(N)} \sum_{\substack{m\leq N/p_1^{a_1}\cdots p_r^{a_r}\\ \gcd(m, N_r)=1}} \sigma_{\text{odd}}(m) \\
& \quad \quad \quad \quad \quad \quad \quad \quad \quad \quad \quad \quad \quad \cdot \sigma_{\text{odd}}(p_1^{a_1}\cdots p_r^{a_r})\mathbf{1}_{G_{N_r}}(p_1^{a_1}\cdots p_r^{a_r}).
\end{align*}
By Lemma \ref{Lemma:NumquatsleN}, 
\begin{align*}
d(S(G)) &\ = \ \lim_{r\rightarrow \infty} \sum_{(a_j)_{j=1}^r \in \mathbb{N}^r} \left( \lim_{N\rightarrow \infty} \frac{\left(\frac{N}{p_1^{a_1}\cdots p_r^{a_r}}\right)^2 \cdot \frac{\pi^2}{24} \frac{\varphi(2N_r)}{N_r} \prod_{j\leq r}(1- p_j^{-2})+ O(N\log N)}{N^2 \cdot \frac{\pi^2}{24} + O(N\log N)} \right) \\
& \quad \quad \quad \quad \quad \quad \quad \quad \quad \quad \quad \quad \quad \cdot \sigma_{\text{odd}}(p_1^{a_1}\cdots p_r^{a_r})\mathbf{1}_{G_{N_r}}(p_1^{a_1}\cdots p_r^{a_r}).
\end{align*}
We see
\begin{align*}
d(S(G)) &\ = \ \lim_{r\rightarrow \infty} \sum_{(a_j)_{j=1}^r \in \mathbb{N}^r} \left(\frac{\varphi(2N_r)}{N_r}\prod_{j\leq r} \frac{1}{p^{2a_j}}(1-p_j^{-2})\right) \cdot \sigma_{\text{odd}}(p_1^{a_1}\cdots p_r^{a_r})\mathbf{1}_{G_{N_r}}(p_1^{a_1}\cdots p_r^{a_r}) \\
&\ = \ \frac{1}{\zeta(2)} \lim_{r\rightarrow \infty} \frac{\varphi(2N_r)}{N_r} \sum_{(a_j)_{j=1}^r \in \mathbb{N}^r}\frac{\sigma_{\text{odd}}(p_1^{a_1}\cdots p_r^{a_r})\cdot \mathbf{1}_{G_{N_r}}(p_1^{a_1}\cdots p_r^{a_r})}{p_1^{2a_1} \cdots p_r^{2a_r}} \\
&\ = \ \frac{1}{\zeta(2)} \lim_{r\rightarrow \infty} \frac{\varphi(2N_r)}{N_r} \prod_{j\leq r} \sum_{k=0}^\infty \frac{\sigma_{\text{odd}}(p_j^k)\cdot \mathbf{1}_{G_{p_j}}(p_j^k)}{p_j^{2k}}\\
&\ = \ \frac{2}{\zeta(2)}\prod_{p \text{ prime}} \left(1 - \frac{1}{p}\right)\sum_{k=0}^\infty \frac{\sigma_{\text{odd}}(p^k)\cdot \mathbf{1}_G(p^k)}{p^{2k}} \ = \ \frac{2}{\zeta(2)} \prod_{p \text{ prime}} \left(1 - \frac{1}{p} \right)\sum_{k\in A} \frac{\sigma_{\text{odd}}(p^k)}{p^{2k}}.
\end{align*}
Using the fact $\sigma_{\text{odd}}(2^k)=1$  and $\sigma_{\text{odd}}(p^k)=\frac{p^{k+1}-1}{p-1}$ for primes $p>2$, we get
\[
d(S(G)) \ = \ \frac{1}{\zeta(2)} \left(\sum_{k\in A} \frac{1}{2^{2k}}\right)\prod_{\substack{p>2\\ \text{prime}}}\left(\sum_{k\in A} \frac{1}{p^k} - \frac{1}{p}\sum_{k\in A} \frac{1}{p^{2k}}\right). \qedhere
\]
\end{proof}

\begin{example}\label{example:squareFreeHurDensity}
If \textit{A} denotes the set of all Hurwitz quaternions with squarefree norm, then by Theorem \ref{thm:eulerProdForQuatDens} and the fact that $n$ is squarefree if and only if $\nu_p(n) \in \{0,1\}$ for all primes $p$, we see 
\begin{align*}
\ d(A) 
&\ = \ \frac{1}{\zeta(2)} \left(1 + \frac{1}{2^{2}}\right)\prod_{\substack{p>2\\ \textup{prime}}}\left(1 + \frac{1}{p} - \frac{1}{p} - \frac{1}{p^{3}}\right) \\
&\ = \ \frac{1}{\zeta(2)\zeta(3)} \frac{1 + 1/2^2}{1- 1/2^3} \ = \ \frac{60}{7\pi^2 \zeta(3)} ~ \approx ~ 0.722484.
\end{align*}
\end{example}

Now, we use Theorem \ref{thm:eulerProdForQuatDens} to compute the density of $S(G_3^*(\mathbb{N}_+))$.

\begin{theorem}\label{Theorem:Rankindensity}
The density of the set $S(G_3^*(\mathbb{N}_+))$ of Hurwitz quaternions with norm in Rankin's greedy set $G_3^*(\mathbb{N}_+)$ is 
$$d(S(G_3^*(\mathbb{N}_+))) \ = \ \frac{f(2^2)}{\zeta(2)}  \prod_{\substack{p >2 \\  \emph{prime}}}\left(f(p)-\frac{f(p^2)}{p}\right)  \ \approx \ 0.782643,$$
where $$f(x) \ = \ \prod_{i=0}^\infty \left(1+\frac{1}{x^{3^i}}\right).$$
\end{theorem}

\begin{proof}
By Theorem \ref{thm:eulerProdForQuatDens} we see 
\begin{align*}
d(S(G_3^*(\mathbb{N}_+))) &\ = \ \frac{1}{\zeta(2)} \left(\sum_{k\in A_3^*(\mathbb{N}_+)} \frac{1}{2^{2k}}\right) \cdot \prod_{\substack{p >2 \\  \text{prime}}}\left(\sum_{k\in A_3^*(\mathbb{N}_+)} \frac{1}{p^k} - \frac{1}{p}\sum_{k\in A_3^*(\mathbb{N}_+)} \frac{1}{p^{2k}}\right) \\
&\ = \ \frac{1}{\zeta(2)}\prod_{j=0}^\infty \left(1 + \frac{1}{2^{2\cdot 3^{j}}}\right) \cdot \prod_{\substack{p >2 \\  \text{prime}}}\left( \prod_{j=0}^\infty \left( 1 + \frac{1}{p^{3^j}}\right) - \frac{1}{p} \prod_{j=0}^\infty \left(1  +\frac{1}{p^{2\cdot 3^j}} \right) \right) \\
&\ = \ \frac{f(2^2)}{\zeta(2)} \prod_{\substack{p >2 \\  \text{prime}}} \left(f(p) - \frac{f(p^2)}{p} \right),
\end{align*}
where in the second line we use the fact that $A_3^*(\mathbb{N}_+)$ is the set of non-negative integers whose ternary expansion do not contain the digit 2 \cite{BG}.

This product converges and is estimated to be $\approx 0.782643$ through Mathematica.
\end{proof}

\section{Bounds on the supremum of upper densities}\label{sec:three}

We  now consider how large a set of Hurwitz quaternions can be while avoiding 3-term geometric progressions.  We make this question precise, as follows.
\begin{question}\label{question:supHurAvoid3GP}
What is the supremum of upper densities of sets of Hurwitz quaternions that avoids 3-term geometric progressions?
\end{question}
Here, the \emph{upper density} of a set of Hurwitz quaternions $A$ in $Q_\text{Hur}$, denoted $\overline{d}(A)$, is defined using the same formula as Equation \eqref{eq:densityDef}, except with $\lim$ replaced by $\limsup$.
We answer this question by obtaining upper and lower bounds for the supremum of upper densities of sets of Hurwitz quaternions avoiding 3-term geometric progressions.

\subsection{Lower bound}\label{lowerbound}

For each $N \in \mathbb{N}_+$, let
$$T_N \ \coloneqq \  \left(\left(\frac{N}{48}, \frac{N}{45}\right] \cup \left(\frac{N}{40}, \frac{N}{36}\right]
\cup \left(\frac{N}{32}, \frac{N}{27}\right] \cup \left(\frac{N}{24}, \frac{N}{12}\right]
\cup \left(\frac{N}{9}, \frac{N}{8}\right] \cup \left(\frac{N}{4}, N \right] \right) \cap \mathbb{N}_+.$$

By Lemma \ref{Lemma:NumquatsleN} we have
\begin{align*}
\abs{S(T_N)} &\ = \    \abs{S((N/4,N])} + \abs{S((N/9,N/8])} + \abs{S((N/24,N/12])} \\
& \quad\quad ~+ \abs{S((N/32,N/27])} + \abs{S((N/40,N/36])} + \abs{S((N/48,N/45])} \\
&\ = \  R_N+ r_N,
\end{align*}
where
$$R_N \ \coloneqq \ \frac{17665627}{18662400}\pi^2N^2 \quad \text{ and } \quad r_N \ \coloneqq \ O(N\log N).$$
Indeed, 
\begin{align*}
\frac{17665627}{18662400} \ =& \ 
\left(1^2- \frac{1}{4^2}\right) + \left(\frac{1}{8^2} - \frac{1}{9^2}\right) + \left(\frac{1}{12^2} - \frac{1}{24^2}\right)\\
& ~ + \left(\frac{1}{27^2} - \frac{1}{32^2}\right) + \left(\frac{1}{36^2} - \frac{1}{40^2}\right) + \left(\frac{1}{45^2} - \frac{1}{48^2}\right).
\end{align*}

Now, write $\mathcal{S}$ to denote $$\mathcal{S} \ \coloneqq \ \bigsqcup_{i=1}^\infty T_{N_i}, \quad \text{ where } \quad N_i = 48^{2^i} \text{ for each $i \in \mathbb{N}$}.$$
McNew \cite[proof of Thm.\ 3.1]{McN} showed that $\mathcal{S}$ is free of geometric progressions with integral ratio.  Thus, $S(\mathcal{S})$ is a set of Hurwitz quaternions that avoids 3-term geometric progressions. We compute the upper density of $S(\mathcal{S})$ in $Q_\text{Hur}$, thereby obtaining a lower bound for the supremum of upper densities of sets of Hurwitz quaternions avoiding $3$-term geometric progressions.

\begin{lemma}\label{Lemma:lowerbound}
The upper density of $S(\mathcal{S})$ in $Q_\textup{Hur}$ is $$\overline{d}(S(\mathcal{S})) \ = \ \frac{17665627}{18662400} \ \approx \  0.946589 .$$
\end{lemma}
\begin{proof}
By the definition of upper density and by Lemma \ref{Lemma:NumquatsleN}, we see that
\begin{align*}
\overline{d}(S(\mathcal{S})) &\ = \ \limsup_{N\rightarrow \infty} \frac{\abs{S(\mathcal{S}) \cap S(N)}}{\abs{S(N)}} \ = \ \lim_{k \to \infty} \frac{\Big|\bigsqcup_{i=1}^{k} S(T_{N_{i}})\Big|}{|S(N_k)| } \notag \\
&\ = \ \lim_{k \to \infty}\frac{R_{N_k}+R_{N_{k-1}}+R_{N_{k-2}}+\cdots+R_{N_1}}{\pi^2N_k^2+O(N_k\textup{log}N_k)}+\lim_{k \to \infty}\frac{r_{N_k}+r_{N_{k-1}}+r_{N_{k-2}}+\cdots+r_{N_1}}{\pi^2N_k^2+O(N_k\textup{log}N_k)}.
\end{align*}
We first argue that the second limit term in the above expression is zero.  Notice that 
$$k\cdot O(1)\ \leq \ r_{N_k}+r_{N_{k-1}}+r_{N_{k-2}} + \cdots+r_{N_1}\ \leq \ k\cdot O\Big(N_k\textup{log}N_k\Big).$$
Since $N_k=48^{2^k}$, we have 
$$\lim_{k\to\infty}\frac{k\cdot O(1)}{\pi^2N_k^2+O(N_k\textup{log}N_k)}\ = \ 0 \quad \text{and} \quad \lim_{k\to\infty}\frac{k\cdot O\Big(N_k\textup{log}N_k\Big)}{\pi^2N_k^2+O(N_k\textup{log}N_k)}\ = \ 0.$$
Hence,
\begin{align*}
\overline{d}(S(\mathcal{S})) &\ = \ \lim_{k \to \infty}\frac{R_{N_k}+R_{N_{k-1}}+R_{N_{k-2}}+\cdots+R_{N_1}}{\pi^2N_k^2+O(N_k\textup{log}N_k)}+0\\
&\ = \ \lim_{k\to\infty}\left(1+\frac{1}{\left(48^{2^{k-1}}\right)^2}+\frac{1}{\left(48^{2^{k-1}+2^{k-2}}\right)^2}+\cdots+\frac{1}{\left(48^{2^{k-1}+2^{k-2}+\cdots+2}\right)^2}\right)\\
& \quad\quad \cdot \lim_{k \to \infty}\frac{R_{N_k}}{\pi^2N_k^2+O(N_k\textup{log}N_k)}.
\end{align*}
Notice that
$$0\ \leq \ \frac{1}{\left(48^{2^{k-1}}\right)^2}+\frac{1}{\left(48^{2^{k-1}+2^{k-2}}\right)^2}+\cdots+\frac{1}{\left(48^{2^{k-1}+2^{k-2}+\cdots+2}\right)^2}\ \leq \ \frac{k-1}{\big(48^{2^{k-1}}\big)^2}\xrightarrow[k\rightarrow \infty]{} 0,$$
so the first limit above converges to 1 and we find that
\[
\overline{d}(S(\mathcal{S})) \ = \ 1 \cdot \lim_{k \to \infty}\frac{R_{N_k}}{\pi^2N_k^2+O(N_k\textup{log}N_k)}\ =\ \frac{17665627}{18662400}.  \qedhere
\]
\end{proof}

\subsection{Upper bound}\label{upperbound} The approach we use to obtain an upper bound is similar to the approach used by McNew \cite{McN} in the classical setting.  Namely, we identify a large collection of \textit{disjoint} 3-term geometric progressions from which at least one term of each progression must be excluded. 
\begin{definition}
Two 3-term geometric progressions of Hurwitz quaternions $(b_1,b_1r_1,b_1r_1^2)$ and $(b_2,b_2r_2,b_2r_2^2)$ are \emph{disjoint} if $$\{b_1,b_1r_1,b_1r_1^2\}\cap \{b_2,b_2r_2,b_2r_2^2\} \ = \ \varnothing.$$
\end{definition}

Fix a Hurwitz quaternion $R$ of norm $2$. For example, we can take $R=1+i$.  For each $n\in\mathbb{N}$ and $N\in\mathbb{N}_+$, write $E_n(N)$ to denote $$E_n(N) \ \coloneqq \ \{ (bR^{3n},bR^{3n+1},bR^{3n+2}) ~\mid ~ 2\nmid \text{Nm}(b) \leq N/2^{3n+2} \}.$$ Note that $E_n(N)$ is empty if $N< 2^{3n+2}$. 
The cardinality of $E_n(N)$ is equal to the number of Hurwitz quaternions with norm less than or equal to $N/2^{3n+2}$ and coprime to $2$.  By Lemma \ref{Lemma:NumquatsleN}, for each fixed $n\in \mathbb{N}$,
$$\abs{E_n(N)} \ = \ 
\pi^2 \left(\frac{N}{2^{3n+2}} \right)^2 \cdot \frac{3}{4} \quad+\quad O(N\log N).$$

\begin{lemma}\label{Lemma:disjoint3TGPs}
The union $$\mathcal{E}(N) \ \coloneqq \ \bigcup_{n=0}^\infty E_n(N)$$  consists of disjoint 3-term geometric progressions contained in $S(N)^3$.  Its cardinality is $$\abs{\mathcal{E}(N)} \ = \ \pi^2 N^2\cdot \frac{1}{21} \ +\  O(N (\log N)^2).$$
\end{lemma}
\begin{proof}
For every Hurwitz quaternion $Q\in Q_\textup{Hur}$, there is a unique way to write $Q$ as $Q = b R^m$ for some $b\in Q_\textup{Hur}$ with $2\nmid \text{Nm}(b)$ and $m\in \mathbb{N}$.  This implies that $\mathcal{E}(N)$ consists of disjoint 3-term geometric progressions. The fact that $\mathcal{E}(N) \subset S(N)^3$ is obvious.

Note that $E_n(N) \cap E_m(N) = \varnothing$ for $n\neq m \in \mathbb{N}$.  Combined with the fact that $E_n(N)$ is empty for $N < 2^{3n+2}$, we see $$\mathcal{E}(N)\ = \ \bigsqcup_{n=0}^{B(N)} E_n(N), \quad  \text{where} \quad B(N)\ \coloneqq \ \left\lceil\tfrac{1}{3}(\log_2(N) - 2)\right\rceil.$$ 
where $\log_2 N$ denotes the logarithm base 2.  Thus, 
\begin{align*}
\abs{\mathcal{E}(N)} &\ = \ \sum_{n=0}^{B(N)}\left( \pi^2 \left(\frac{N}{2^{3n+2}} \right)^2 \cdot \frac{3}{4} \ +\ O(N \log N) \right) \\
&\ = \ \pi^2 N^2 \cdot \frac{3}{4^3}\sum_{n=0}^{B(N)} \left(\frac{1}{4^{3}}\right)^n \ +\ O(N(\log N)^2) \ = \ \pi^2N^2 \cdot \frac{1}{21} + O(N (\log N)^2).
\qedhere\end{align*}
\end{proof}

\begin{lemma}\label{Lemma:upperbound}
Let $A\subset Q_\textup{Hur}$ be any set of Hurwitz quaternions avoiding 3-term geometric progressions.  Then $$\overline{d}(A) \ \leq \ \frac{20}{21} \ \approx \ 0.952381.$$
\end{lemma}

\begin{proof}
Observe $\abs{A\cap S(N)} = \abs{S(N)} - \abs{S(N) \setminus A}$.  Since, $\mathcal{E}(N)$ consists of 3-term geometric progressions in $S(N)$, at least one term from each element of $\mathcal{E}(N)$ must be contained in $S(N) \setminus A$.  Since all of the 3-term geometric progressions in $\mathcal{E}(N)$ are disjoint, the cardinality of $S(N)\setminus A$ is at least the cardinality of $\mathcal{E}(N)$.  By Lemma \ref{Lemma:disjoint3TGPs}, $$\abs{S(N)\setminus A} \ \geq \ \pi^2 N^2\cdot \frac{1}{21} \ +\  O(N (\log N)^2).$$

Thus,
\begin{align*}
\overline{d}(A)&\ = \ \limsup_{N\rightarrow\infty}\frac{\abs{A\cap S(N)}}{\abs{S(N)}} \ = \ \limsup_{N\rightarrow\infty} \frac{\abs{S(N)} - \abs{S(N) \setminus A}}{\abs{S(N)}} \\
&\ \leq \ 1 - \liminf_{n\rightarrow \infty} \frac{\pi^2 N^2 \cdot \frac{1}{21} \ +\ O(N(\log N)^2)}{\pi^2 N^2 \ +\ O(N\log N)} \ = \ \frac{20}{21}.
\qedhere\end{align*}
\end{proof}
We have finally arrived at our answer to Question \ref{question:supHurAvoid3GP}.
\begin{theorem}\label{mainthm-quatupperlowerbd}
The supremum $\overline{m}_{\textup{Hur}}$ of upper densities of sets of Hurwitz quaternions that avoid 3-term geometric progressions satisfies the following bounds:
$$ 0.946589 \ \approx \ \frac{17665627}{18662400} \ \le\ \overline{m}_{\textup{Hur}}\ \le\ \frac{20}{21} \ \approx \ 0.952381.$$
\end{theorem}
\begin{proof}
Combine Lemma \ref{Lemma:lowerbound} and Lemma \ref{Lemma:upperbound}.
\end{proof}

\section{The quaternion greedy set} \label{Section: Greedy Set}
In Section \ref{Section: Rankin}, we computed an exact formula for the density of the set of Hurwitz quaternions $S(G_3^*(\mathbb{N}_+))$ with norm in Rankin's greedy set $G_3^*(\mathbb{N}_+)$.  It is natural to consider a different set of Hurwitz quaternions avoiding 3-term geometric progressions that is constructed by a greedy algorithm in a different way.  Namely, consider the set of Hurwitz quaternions formed by including Hurwitz quaternions of increasing norms so long as they do not form a geometric progression with elements of smaller norms already included in the set.   This process begins with including all the units, or Hurwitz quaternions of norm 1, and then considers progressively larger norms.  The set constructed in this way will be well-defined because including a particular quaternion of a given norm, $n$, will not create a geometric progression with any other quaternions of norm $n$ since unit ratios are not allowed.
We call this set the \emph{quaternion greedy set} and denote it as $G_3^*(Q_\textup{Hur}\setminus\{0\})$.  The following observation is an immediate consequence of this definition.
\begin{proposition}\label{prop:suffCondQuadGreedy}
Any Hurwitz quaternion which is not the third term of a 3-term geometric progression is contained in the quaternion greedy set $G_3^*(Q_\textup{Hur}\setminus\{0\})$.
\end{proposition}
\begin{question}\label{Question:QuatGreedyDensity}
What is the (upper/lower) density of $G_3^*(Q_\textup{Hur}\setminus\{0\})$?
\end{question}

We are not able to answer Question \ref{Question:QuatGreedyDensity} completely.  However, the following proposition gives a lower bound for the density of $G_3^*(\mathbb{N}_+)$.

\begin{proposition}\label{Prop:LowerBoundQuatGreedy}
The quaternion greedy set $G_3^*(Q_\textup{Hur}\setminus\{0\})$ contains the set of Hurwitz quaternions whose norm is squarefree.  Thus, the lower density of $G_3^*(Q_\textup{Hur}\setminus\{0\})$ in $Q_\textup{Hur}$ is at least $$\underline{d}(G_3^*(Q_\textup{Hur}\setminus\{0\})) \ \geq \ \frac{60}{7\pi^2 \zeta(3)} \ \approx \ 0.722484.$$
\end{proposition}

\begin{proof}
If $Q$ is a Hurwitz quaternion whose norm is squarefree, then $Q$ cannot be the third term in a 3-term geometric progression, and thus by Proposition \ref{prop:suffCondQuadGreedy} must be included in the quaternion greedy set.  The lower bound for the density follows from Example \ref{example:squareFreeHurDensity}.
\end{proof}

In the remainder of this section, we attempt to convince the reader that Question \ref{Question:QuatGreedyDensity} is hard.  One may try to tackle Question \ref{Question:QuatGreedyDensity} via the same method used to compute the density of $S(G_3^*(\mathbb{N}_+))$.  However, step zero of this method would be to realize $G_3^*(Q_\textup{Hur}\setminus\{0\})$ as the set of Hurwitz quaternions whose norm is contained in a given set $G \subset \mathbb{N}_+$.  This is not possible, as shown by the following proposition.

\begin{proposition}\label{observation}
The quaternion greedy set is not equal to $S(G)$ for any $G\subset \mathbb{N}_+$.
\end{proposition}

We postpone the proof of Proposition \ref{observation} to the end of this section.

Another approach to Question \ref{Question:QuatGreedyDensity} is to try and relate the quaternion greedy set $G_3^*(Q_\textup{Hur}\setminus\{0\})$ to $S(G_3^*(\mathbb{N}_+))$.  One might hope that $G_3^*(Q_\textup{Hur}\setminus\{0\})$ would contain or be contained within $S(G_3^*(\mathbb{N}_+))$ (with possibly finitely many exceptions).  If this were true, then one would be able to bound the density of $G_3^*(Q_\textup{Hur}\setminus\{0\})$ in terms of the already computed density of $S(G_3^*(\mathbb{N}_+))$. However, this is false, as shown by the following proposition.

\begin{proposition}\label{Prop:InfinitelyManyComplements}
There are infinitely many Hurwitz quaternions in both complement sets $G_3^*(Q_\textup{Hur}\setminus\{0\}) \setminus S(G_3^*(\mathbb{N}_+))$ and $S(G_3^*(\mathbb{N}_+)) \setminus G_3^*(Q_\textup{Hur}\setminus\{0\})$.
\end{proposition}

We postpone the proof of Proposition \ref{Prop:InfinitelyManyComplements} until the end of this section.  The proof will provide infinitely many such Hurwitz quaternions explicitly.

A key step in establishing both Proposition \ref{observation} and Proposition \ref{Prop:InfinitelyManyComplements} is realizing that there exist Hurwitz quaternions of norm $r^2 > 0$ that cannot be written as the square of a Hurwitz quaternion of norm $r$ multiplied by a unit on the left. 

\begin{proposition}\label{Prop:generalproposition}
Suppose $r\in\mathbb{N}_+$ is a positive integer such that either
\begin{itemize}
    \item $r$ cannot be represented as the sum of three integer squares, or
    \item $r$ is divisible by an odd integer greater than 24.
\end{itemize} Then there exists a Hurwitz quaternion $Q$ of norm $\textup{Nm}(Q)=r^2$ such that $Q$ cannot be written in the form $Q=UR^2$ for any unit $U\in Q_{\textup{Hur}}$ and any $R\in Q_{\textup{Hur}}$ of norm $\textup{Nm}(R)=r$.  In the first case, $Q=r$ is such a Hurwitz quaternion.
\end{proposition}

\begin{proof}[Proof when $r$ is not a sum of 3 squares]
Let $r\in\mathbb{N}_+$ be such that $r$ cannot be represented as the sum of three integer squares.  Suppose for contradiction that $r=UR^2$ for some unit $U\in Q_{\text{Hur}}$ and $R\in Q_{\text{Hur}}$ of  $\textup{Nm}(R)=r$.  Then $U^{-1}r=R^2$. Suppose that $U^{-1}=a+bi+cj+dk$ and $R=w+xi+yj+zk$.  Then, solving the equation $U^{-1}r=R^2$, we find that $a=\frac{w^2-x^2-y^2-z^2}{r}$, $b=\frac{2wx}{r}$, $c=\frac{2wy}{r}$, $d=\frac{2wz}{r}$.  It is a consequence of Lagrange's three-square theorem that the set of integers that cannot be represented as a sum of three integer squares is closed under multiplication by $4$.  Since $\textup{Nm}(R)=r$ cannot be represented as the sum of three integer squares, we must have that $w, x, y, z$ are all nonzero. Since $w, x, y, z$ are nonzero, then $b, c, d$ are nonzero; moreover, as $U^{-1}=a+bi+cj+dk$ is a unit in $Q_{\text{Hur}}$, then $a$ is also nonzero, thus $a, b, c, d\in \{\pm\frac{1}{2}\}$. Then, $x=\pm\frac{r}{4w}$, $y=\pm\frac{r}{4w}$, and $z=\pm\frac{r}{4w}$. Substituting into $\frac{w^2-x^2-y^2-z^2}{r}=\pm\frac{1}{2}$ and solving the equation, we find that $w=\pm\frac{\sqrt{r}}{2}$ or $w=\pm\frac{\sqrt{3r}}{2}$. We know that $r$ cannot be represented as the sum of three integer squares; hence, $r$ is not an integer square and $w\neq\pm\frac{\sqrt{r}}{2}$. We also know $w=\pm\frac{\sqrt{3r}}{2}\in\mathbb{Z}$ or $\mathbb{Z}+\frac{1}{2}$ if and only if $r=3l^2$ for some $l\in\mathbb{Z}$. However, this implies that $r$ can be represented as the sum of three integer squares, which is a contradiction.
\end{proof}
\begin{proof}[Proof when $r$ has an odd divisor $> 24$]
By Lemma \ref{Lemma:NumquatsfixedN}, the number of Hurwitz quaternions $Q$, which can be written as $Q= UR^2$ where $U \in Q_\textup{Hur}$ is a unit and $R\in Q_\textup{Hur}$, has norm $\textup{Nm}(R) = r$ is $24^2 \cdot \sigma_\text{odd}(r)$.  Thus, to show there exists Hurwitz quaternions $Q$ of norm $\textup{Nm}(Q) = r^2$, which cannot be written in this form, it suffices to show $24 \cdot \sigma_\text{odd}(r) \ < \ \sigma_\text{odd}(r^2)$.

Let $D$ be the greatest odd divisor of $r$.  By hypothesis, $D>24$. Thus, we have $$24 \cdot \sigma_\textup{odd}(r) \ < \ \sum_{2 \nmid d \mid r} D\cdot d \ \le \ \sum_{2 \nmid d \mid r^2} d \ = \ \sigma_\textup{odd}(r^2),$$
as required.
\end{proof}
\begin{proposition}\label{Prop:Complements}
Let $p$ be a rational prime number. Suppose $Q$ is a Hurwitz quaternion of norm $\textup{Nm}(Q) = p^2$ such that $Q$ cannot be written in the form $Q = UR^2$ for any unit $U \in Q_\textup{Hur}$ and any $R\in Q_\textup{Hur}$ of norm $\textup{Nm}(R) = p$.  Then $Q\in G_3^*(Q_\textup{Hur}\setminus\{0\})\setminus S(G_3^*(\mathbb{N}_+))$ and $Q^2 \in S(G_3^*(\mathbb{N}_+)) \setminus G_3^*(Q_\textup{Hur}\setminus\{0\})$.
\end{proposition}
\begin{proof}
In light of Proposition \ref{prop:suffCondQuadGreedy}, in order to prove that $Q \in G_3^*(Q_\textup{Hur}\setminus\{0\})$ it suffices to prove that $Q$ is not the third term of any 3-term geometric progression.  Suppose for contradiction that $Q= br^2$ for some $b,r \in Q_\textup{Hur} \setminus \{0\}$ with $r$ a non-unit.  Then we would have $p^2 = \textup{Nm}(Q) = \textup{Nm}(b)\cdot \textup{Nm}(r)^2$.  Since $p$ is prime, $\textup{Nm}(b),\textup{Nm}(r)\in \mathbb{N}_+$ and $\textup{Nm}(r) > 1$, we see $\textup{Nm}(r) = p$ and $\textup{Nm}(b) = 1$.  Thus, $b$ is a unit, and $r$ has norm $p$, contradicting the supposition that $Q$ cannot be written in the form $Q=UR^2$.  Thus, $Q\in  G_3^*(Q_\textup{Hur}\setminus\{0\})$.  

To show that $Q\in G_3^*(Q_\textup{Hur}\setminus\{0\})\setminus S(G_3^*(\mathbb{N}_+))$, it remains to prove $Q\not\in S(G_3^*(\mathbb{N}_+))$. Since $\text{Nm}(Q) = p^2$, and $2\not\in A_3^*(\mathbb{N}_+)$, we know $Q\not\in S(G_3^*(\mathbb{N}_+))$.  

The norm of $Q^2$ is $\textup{Nm}(Q^2) = p^4$ and $4 \in A_3^*(\mathbb{N}_+)$, so $Q^2 \in S(G_3^*(\mathbb{N}_+))$.  We already know $Q \in G_3^*(Q_\textup{Hur}\setminus\{0\})$. Since both the first and second terms of the 3-term geometric progression $(1,Q,Q^2)$ are in $G_3^*(Q_\textup{Hur}\setminus\{0\})$, we see that $Q^2 \not\in G_3^*(Q_\textup{Hur}\setminus\{0\})$.  Thus, $Q^2\in S(G_3^*(\mathbb{N}_+)) \setminus G_3^*(Q_\textup{Hur}\setminus\{0\}) $.
\end{proof}
Armed with Proposition \ref{Prop:generalproposition} and Proposition \ref{Prop:Complements}, we can now prove Proposition \ref{observation} and Proposition \ref{Prop:InfinitelyManyComplements}.

\begin{proof}[Proof of Proposition \ref{observation}]    
Since $7$ is not a sum of three squares, we know by Proposition \ref{Prop:generalproposition} and Proposition \ref{Prop:Complements} that $7 \in G_3^*(Q_\textup{Hur}\setminus\{0\})$.  Note $Q=2+i+j+k$ is a Hurwitz quaternion of norm $7$.  By Proposition \ref{Prop:LowerBoundQuatGreedy}, $Q \in G_3^*(Q_\textup{Hur}\setminus\{0\})$.  Since both $1$ and $Q$ are in $G_3^*(Q_\textup{Hur}\setminus\{0\})$, the third term $Q^2$ of the 3-term geometric progression $(1,Q,Q^2)$ cannot be in $G_3^*(Q_\textup{Hur}\setminus\{0\})$.  Thus, $G_3^*(Q_\textup{Hur}\setminus\{0\})$ contains some Hurwitz quaternions of norm $7^2$, but not all.  Therefore, $G_3^*(Q_\textup{Hur}\setminus\{0\})$ cannot be $S(G)$ for some $G \subset \mathbb{N}_+$.
\end{proof}

\begin{proof}[Proof of Proposition \ref{Prop:InfinitelyManyComplements}]
This follows from Proposition \ref{Prop:generalproposition} and Proposition \ref{Prop:Complements}.  Note that as $p$ ranges over the infinitely many primes that cannot be represented as the sum of three integer squares, Proposition \ref{Prop:generalproposition} and Proposition \ref{Prop:Complements} imply that $p\in G_3^*(Q_\textup{Hur}\setminus\{0\}) \setminus S(G_3^*(\mathbb{N}_+))$ and $p^2 \in S(G_3^*(\mathbb{N}_+)) \setminus  G_3^*(Q_\textup{Hur}\setminus\{0\})$.  This provides an explicit list of infinitely many examples in both complements.
\end{proof}

In conclusion, there is no clear relationship between $G_3^*(Q_\textup{Hur}\setminus\{0\})$ and $S(G_3^*(\mathbb{N}_+))$.  One may try to systematically describe the inclusions and exclusions necessary to transform $S(G_3^*(\mathbb{{Z}}))$ to $G_3^*(Q_\textup{Hur}\setminus\{0\})$, but this appears to be hard to predict or keep track of.  Consequently, we do not know whether the density of $G_3^*(Q_\textup{Hur}\setminus\{0\})$ in $Q_\textup{Hur}$ is greater than or less than the density of $S(G_3^*(\mathbb{N}_+))$ in $Q_\textup{Hur}$.  We also cannot make any computational estimate of the density of the quaternion greedy set as we do not even have a quick way of testing whether a given Hurwitz quaternion is included in $G_3^*(Q_\textup{Hur}\setminus\{0\})$.

\section{Future work}\label{section:FutureWork}

As discussed in Section \ref{Section: Greedy Set}, the nature of the quaternion greedy set $G_3^*(Q_\textup{Hur}\setminus\{0\})$ is quite mysterious.  Fundamentally, it is the noncommutativity of the Hurwitz quaternions which makes the quaternion greedy set so difficult to analyze and understand.  We record here several questions about $G_3^*(Q_\textup{Hur}\setminus\{0\})$ that remain open.

\begin{question}\label{question:naive1}
Does $G_3^*(Q_\textup{Hur}\setminus\{0\})$ have an asymptotic density?
\end{question}

\begin{question}
Is the (lower/upper) density of $G_3^*(Q_\textup{Hur}\setminus\{0\})$ greater than $d(S(G_3^*(\mathbb{N}_+)))$?
\end{question}

\begin{question}\label{question:naive4}
Is there an alternative characterization of the quaternion greedy set?
Is there an efficient algorithm that tests membership in $G_3^*(Q_\textup{Hur}\setminus\{0\})$?
\end{question}

Next, we introduce a general framework in which questions about avoiding geometric progressions can be studied.  Let $B$ be a cancellative, normed monoid.  This means $B$ is a monoid which satisfies the cancellation property, and $B$ comes equipped with a monoid homomorphism $N: B \rightarrow \mathbb{N}_+$ which we call the norm.  Suppose that for each $N\in \mathbb{N}_+$, the set $S_B(N)$ of elements in $B$ with norm $\leq N$ is finite.  We define the \textit{asymptotic density of a set} $A \subset B$ in $B$ to be $$d_B(A) \ \coloneqq \ \lim_{N\rightarrow \infty} \frac{\abs{A \cap S_B(N)}}{\abs{S_B(N)}}$$ if it exists.  The upper density $\overline{d}_B(A)$ is defined similarly, except $\lim$ is replaced with $\limsup$.

Let $U(B)$ be the universal enveloping group.  This is the group with generators given by the elements of $B$ and relations given by the full multiplication table of $B$.  Since $B$ is cancellative, the natural monoid homomorphism $B\rightarrow U(B)$ is injective, and we regard $B$ as a subset of $U(B)$ via this map. The norm $N: B \rightarrow \mathbb{N}_+$ naturally extends to a group homomorphism $U(N) : U(B) \rightarrow \mathbb{Q}_+$, which we also call the norm.

Let $R\subset U(B)$ be a subset.  A \emph{$k$-term geometric progression in $B$ with ratio in $R$} is a $k$-tuple $(b_1,b_2,\cdots, b_k)\in B^k$, where 
$$b_1^{-1}b_2 \ = \ b_2^{-1}b_3 \ = \ \ldots \ = \ b_{k-1}^{-1}b_k \ \in R.$$
A set $S\subset B$ is said to \emph{avoid $k$-term geometric progressions with ratio in $R$} if $S^k$ does not contain any such $k$-tuples.  Now we define the constants
\begin{align*}
m_{k}(B,R) \ &\coloneqq \ \sup\left\{ d_B(S) \,\middle|\, 
\begin{array}{c}
\text{$S\subset B$ avoids $k$-term geometric progressions}\\
\text{ with ratio in $R$, and $d_B(S)$ exists}
\end{array}
\right\} \text { and } \\
\overline{m}_{k}(B,R) \ &\coloneqq \ \sup \left\{\overline{d}_B(S) \mid S\subset B \text{ avoids $k$-term geometric progressions with ratio in $R$}\right\}.
\end{align*}

It is clear from the definitions that $m_k(B,R) \leq \overline{m}_k(B,R)$.  Moreover, both $m_k(B,R)$ and $\overline{m}_k(B,R)$ are increasing with respect to $k$ and decreasing with respect to $R$.

\begin{question}\label{question:genQues}
What is the value of $m_{k}(B,R)$ and $\overline{m}_{k}(B,R)$ for different choices of length $k$, normed monoid $B$  and forbidden ratio set $R$?
\end{question}

In this paper, we obtained upper and lower bounds for $m_{k}(B,R)$ and $\overline{m}_{k}(B,R)$  when $k=3$, $B=Q_\textup{Hur}\setminus\{0\}$ and $R = Q_\textup{Hur}\setminus S(1)$.  Indeed, we have
\begin{align*}
    0.782643 \ \leq \ m_\textup{Hur} &\ \coloneqq \  m_3(Q_\textup{Hur}\setminus\{0\},Q_\textup{Hur}\setminus S(1)), \text{ and}  \tag{Theorem \ref{Theorem:Rankindensity}} \\
    0.946589 \ \leq \ \overline{m}_\textup{Hur} & \ \coloneqq \  \overline{m}_3(Q_\textup{Hur}\setminus\{0\},Q_\textup{Hur}\setminus S(1)) \ \leq \ 0.952381. \tag{Theorem \ref{mainthm-quatupperlowerbd}}
\end{align*}

For any $G\subset \mathbb{N}$, write $S_B(G)$ to denote the set of elements of $B$ with norm in $G$.  Since $G_k^*(\mathbb{N}_+)$ avoids $k$-term geometric progressions with rational ratio, we get that $S_B(G_k^*(\mathbb{N}_+))$ avoids $k$-term geometric progressions with ratio in $R$ for any forbidden ratio set $R \subset U(B)$.  Thus, if $d_B(S_B(G_k^*(\mathbb{N}_+)))$ exists, it serves as a lower bound for $m_k(B,R)$ for all possible forbidden ratio sets $R$.
\begin{question}\label{question:genDensityRankinGreedy}
Does $S_B(G_k^*(\mathbb{N}_+))$ have an asymptotic density in $B$?  If yes, compute it.
\end{question}

Theorem \ref{Theorem:Rankindensity} answers Question \ref{question:genDensityRankinGreedy} when $k=3$ and $B=Q_\textup{Hur}\setminus\{0\}$.

Now we generalize the definition of the quaternion greedy set.  Suppose $R$ is a subset of $U(B)$ that avoids elements of norm $1$.  We define a \textit{greedy set} of elements of $B$ avoiding $k$-term geometric progressions with ratio in $R$ as follows.  Start with the set $S_B(1)$ of elements of $B$ of norm $1$ and adjoin to it elements of $B$ of increasing norm so long as the enlarged set continues to avoid $k$-term geometric progressions with ratio in $R$.  This is a well-defined procedure because $R$ does not contain any elements of norm $1$.  Denote the resulting greedy set as $G_{k}^*(B,R)$. Note $G_3^*(\mathbb{N}_+) = G_3^*(\mathbb{N}_+,\mathbb{N}_{>1}) = G_3^*(\mathbb{N}_+, \mathbb{Q}_{>1})$, as observed by Brown and Gordon \cite{BG} and $G_3^*(Q_\textup{Hur}\setminus\{0\}) = G_3^*(Q_\textup{Hur}\setminus\{0\}, Q_\textup{Hur}\setminus S(1))$.

\begin{question}
Repeat Questions \ref{question:naive1} - \ref{question:naive4} with $G_3^*(Q_\textup{Hur}\setminus\{0\})$ and $S(G_k^*(\mathbb{N}_+))$ replaced by $G_k^*(B,R)$ and $S_B(G_k^*(\mathbb{N}_+))$ respectively.  
\end{question}

It would be interesting to study all the above questions when $B = \mathcal{O} \setminus \{0\}$ where $\mathcal{O}$ is a maximal order in the octonions $\mathbb{O}$, as this is a nonassociative setting.  Even in the Hamilton quaternion setting $B = Q_\textup{Hur}\setminus \{0\}$, it would be interesting to vary the forbidden ratio set $R$ from $Q_\textup{Hur} \setminus S(1)$ to $\mathbb{N}_{>1}^\times$, $\mathbb{Q}_{>1}^\times$ or $\{Q \in \mathbb{H}^\times \mid \textup{Nm}(Q) > 1\}$.

\section{Acknowledgments}
This research was conducted as part of the 2015 SMALL REU program and the 2022 Polymath Jr.\ program. The authors were partially supported by NSF grants DMS 1347804, DMS 1265673, DMS 1561945, and DMS 2218374 and by Williams College.  
The authors thank Emma Dinowitz for helpful discussions regarding the proof of Theorem \ref{thm:eulerProdForQuatDens} and the suggestions of the anonymous referee which greatly improved the exposition of the paper.

\bigskip
\hrule
\bigskip

\noindent 2020 \emph{Mathematics Subject Classification}: Primary 11B05; Secondary 05D10, 11Y60.

\noindent \emph{Keywords:} Ramsey theory, greedy algorithm, non-commutative group, quaternion, Hurwitz quaternion, density, asymptotic density, geometric progression, geometric-progression-free

\bigskip
\hrule
\bigskip

\noindent 
(Concerned with sequences \seqnum{A000452}, \seqnum{A004011} and \seqnum{A005836}.)

\end{document}